\documentclass[reqno,a4paper]{amsart}

\usepackage[english]{babel}
\usepackage{amsfonts}
\usepackage{amsmath}
\usepackage{amssymb}
\usepackage{amsthm}
\usepackage{latexsym}
\usepackage{amstext}
\usepackage{enumerate}
\usepackage{graphics}
\usepackage[all]{xy}

\theoremstyle{plain}
\newtheorem{theorem}{Theorem}[section]

\newtheorem{proposition}[theorem]{Proposition}

\theoremstyle{definition}
\newtheorem{definition}[theorem]{Definition}

\newtheorem{remark}[theorem]{Remark}

\numberwithin{equation}{section}

\DeclareMathOperator{\Ric}{\mathrm{Ric}}

\begin{document}

\title[The $\eta$-Einstein condition on indefinite $\mathcal{S}$-manifolds]
{The $\eta$-Einstein condition on indefinite $\mathcal{S}$-manifolds}
\author{Letizia Brunetti}

\newcommand{\acr}{\newline\indent}

\address{Department of Mathematics\acr
         University of Bari ``Aldo Moro''\acr
         Via E. Orabona, 4\acr
         70125 -- Bari\acr
         ITALY}

\email{brunetti@dm.uniba.it}%


\begin{abstract}
An $\eta$-Einstein condition is introduced in the context of indefinite $g.f.f$-manifolds, and a few Schur-type lemmas for indefinite $\mathcal{S}$-manifolds are provided.
\end{abstract}

\maketitle

\section{Introduction and Preliminaries}
The notion of $f$-structure on a $(2n+s)$-dimensional manifold $M$, i.e.~a $(1,1)$-type tensor field $\varphi$ on $M$ of constant rank $2n$ such that $\varphi^3+\varphi=0$, was firstly introduced in 1963 by K.~Yano (\cite{Y01}) as a generalization of both (almost) contact (for $s=1$) and (almost) complex structures (for $s=0$). During the subsequent years, this notion has been furtherly developed by several authors (\cite{Bl01}, \cite{BLY}, \cite{GY}, \cite{GY02}, \cite{IS01}, \cite{N01}, \cite{N02}). Among them, H.~Nakagawa in \cite{N01} and \cite{N02} introduced the notion of framed $f$-manifold, later developed and studied by S.I.~Goldberg and K.~Yano (\cite{GY}, \cite{GY02}) and others with the denomination of globally framed $f$-manifolds.

A manifold $M$ is said to be a \emph{globally framed $f$-manifold} (briefly \emph{$g.f.f$-manifold}) if it carries a globally framed $f$-structure, that is an $f$-structure $\varphi$ such that the subbundle $\ker(\varphi)$ is parallelizable. If $\mathrm{rank}(\ker(\varphi))=s\geqslant1$ the existence of a $g.f.f$-structure on $M$ is equivalent to the existence of $s$ linearly independent global vector fields $\xi_\alpha$ and $1$-forms $\eta^\alpha$, $\alpha\in\{1,\ldots,s\}$, such that
\begin{equation}\label{eq:02}
\varphi^2=-I+\eta^\alpha\otimes\xi_\alpha\qquad\text{and}\qquad\eta^\alpha(\xi_\beta)=\delta^\alpha_\beta,
\end{equation}
where $I$ is the identity mapping. We point out that this kind of structure is also known as ``$f$-structure with complemented frames'' (\cite{Bl01}, \cite{CFF}), or ``almost $r$-contact structure'' (\cite{Van}).

From (\ref{eq:02}) it follows that $\varphi\xi_\alpha=0$ and $\eta^\alpha\circ\varphi=0$, for any $\alpha\in\{1,\ldots,s\}$. Moreover, $TM=\mathrm{Im}(\varphi)\oplus\mathrm{span}(\xi_1,\ldots,\xi_s)$, where $\mathrm{Im}(\varphi)$ is a distribution on $M$ of even rank $r=2n$ on which $\varphi$ acts as an almost complex tensor field, so that $\dim(M)=2n+s$. Each $\xi_\alpha$ is said to be a \emph{characteristic vector field} of the structure. A $g.f.f$-manifold $(M,\varphi,\xi_\alpha,\eta^\alpha)$ is called \emph{normal} if the $(1,2)$-type tensor field $N=[\varphi,\varphi]+2d\eta^\alpha\otimes\xi_\alpha$ vanishes identically (\cite{IS01}).

Globally framed $f$-structures can always be considered together with an associated Riemannian metric (\cite{Bl01}, \cite{Y01}), while for general indefinite metrics some restrictions on the signature have to be observed. Such restrictions disappear in the case of $g.f.f$-manifolds endowed with Lorentzian metrics (see p.\ 214 of \cite{DB}). More recently, a study of a particular class of $g.f.f$-manifolds endowed with an indefinite metric has been carried out in (\cite{LP}). Following \cite{DB, DS, LP}, we say that an indefinite metric $g$ on a $g.f.f$-manifold $(M,\varphi,\xi_\alpha,\eta^\alpha)$ is \emph{compatible} with the $g.f.f$-structure $(\varphi,\xi_\alpha,\eta^\alpha)$ if
\begin{equation}\label{eq:03}
	g(\varphi X,\varphi Y) = g(X,Y)-\sum_{\alpha=1}^s\varepsilon_\alpha\eta^\alpha(X)\eta^\alpha(Y),
\end{equation}
for all $X,Y\in\Gamma(TM)$, where $\varepsilon_\alpha=g(\xi_\alpha,\xi_\alpha)=\pm1$. Then, the manifold $M$ is said to be an \emph{indefinite (metric) $g.f.f$-manifold} with structure $(\varphi,\xi_\alpha,\eta^\alpha,g)$. From (\ref{eq:03}) we easily get
\begin{equation}\label{eq:11}
g(X,\xi_\alpha)=\varepsilon_\alpha\eta^\alpha(X) \qquad \text{and} \qquad g(X,\varphi Y)=-g(\varphi X,Y),
\end{equation}
for any $X,Y\in\Gamma(TM)$ and any $\alpha\in\{1,\ldots,s\}$. Furthermore, $\mathrm{Im}(\varphi)$ is orthogonal to $\mathrm{span}(\xi_1,\ldots,\xi_s)$, and since $g(\varphi X,\varphi Y)=g(X,Y)$, for any $X,Y\in\mathrm{Im}(\varphi)$, then the signature of $g$ on $\mathrm{Im}(\varphi)$ is $(2p,2q)$, with $2p+2q=2n$. In \cite{DB} it is proved that there always exists a Lorentzian metric $g$ associated with a $g.f.f$-manifold, and in this case exactly one of the characteristic vector fields has to be unit timelike and the restriction of $g$ to $\mathrm{Im}(\varphi)$ has Riemannian signature.

The $2$-form $\Phi$ on $M$ defined by $\Phi(X,Y)=g(X,\varphi Y)$ is called the \emph{fundamental $2$-form} of the indefinite $g.f.f$-manifold. If $\Phi=d\eta^\alpha$, for any $\alpha\in\{1,\ldots,s\}$, the manifold $(M,\varphi,\xi_\alpha,\eta^\alpha,g)$ is said to be an \emph{indefinite almost $\mathcal{S}$-manifold}. Finally, a normal indefinite almost $\mathcal{S}$-manifold is, by definition, an \emph{indefinite $\mathcal{S}$-manifold}. As proved in \cite{LP}, in an indefinite $\mathcal{S}$-manifold the covariant derivative of $\varphi$ satisfies the identity
\begin{equation}\label{eq:04}
(\nabla_X\varphi)Y=g(\varphi X,\varphi Y)\bar{\xi}+\bar{\eta}(Y)\varphi^2X,
\end{equation}
where $\bar{\xi}=\sum_{\alpha=1}^s\xi_\alpha$ and $\bar{\eta}=\sum_{\alpha=1}^s\varepsilon_\alpha\eta^\alpha$, from which it easily follows that, for any $\alpha,\beta\in\{1,\ldots,s\}$,
\begin{equation}\label{eq:12}
\nabla_X\xi_\alpha=-\varepsilon_\alpha\varphi X \qquad \text{and} \qquad \nabla_{\xi_\alpha}\xi_\beta=0,
\end{equation}
as well as that each $\xi_\alpha$ is a Killing vector field. In particular, for $s=1$ one finds again the notion of indefinite Sasakian manifold (\cite{T}).

In the Riemannian setting, the notion of $\mathcal{S}$-manifold, together with other remarkable classes of $g.f.f$-manifolds, appears in \cite{Bl01}. It has been developed by several authors and for further properties we refer the reader to \cite{Bl01}, \cite{BLY}, \cite{CFF} and \cite{DIP}, where the notion of almost $\mathcal{S}$-manifold is introduced. The generalization of this notion to the semi-Riemannian setting is given in \cite{LP}.

The main purpose of this short note is to extend the notion of $\eta$-Einstein $g.f.f$-structure to the semi-Riemannian setting, by a suitable generalization of the definition given in \cite{KT}. We provide it in Section $2$, pointing out the main differences between our definition and that contained in \cite{KT}, and give a first Schur-type lemma. Based on it, in Section $3$, we prove a second Schur-type result for indefinite $\mathcal{S}$-manifolds with pointwise constant $\varphi$-sectional curvature and its suitable consequence. In a forthcoming paper, we are going to develop and apply the results contained here to the study of the $\varphi$-null Osserman condition on Lorentzian $\mathcal{S}$-manifolds.

In what follows, all manifolds, tensor fields and maps are assumed to be smooth. Moreover, all manifolds are supposed to be connected and, according to \cite{KN}, for the curvature tensors of a semi-Riemannian manifold $(M,g)$ we put
\[
R(X,Y,Z,W)=g(R(Z,W)Y,X)=g(([\nabla_Z,\nabla_W]-\nabla_{[Z,W]})Y,X),
\]
for any vector fields $X,Y,Z,W$ on $M$. Finally, for any $p\in M$ and any linearly independent vectors $x,y\in T_pM$ spanning a non-degenerate plane $\pi=\mathrm{span}(x,y)$, that is $\Delta(\pi)=g_p(x,x)g_p(y,y)-g_p(x,y)^2\neq0$, the sectional curvature of $(M,g)$ at $p$ with respect to $\pi$ is, by definition, the real number
\[
k_p(\pi)=k_p(x,y)=\frac{R_p(x,y,x,y)}{\Delta(\pi)}.
\]

\section{The $\eta$-Einstein condition for indefinite $g.f.f$-manifolds.}
Let us now state some preliminary properties of the curvature tensor field of an indefinite $\mathcal{S}$-manifold.

\begin{proposition}
Let $(M,\varphi,\xi^\alpha,\eta_\alpha,g)$, $1\leqslant\alpha\leqslant s$, be a $(2n+s)$-dimensional indefinite $\mathcal{S}$-manifold. The following identities hold, for any $X,Y,Z\in\Gamma(TM)$, any $U,V\in\mathrm{span}(\xi_1,\ldots,\xi_s)$ and any $\alpha,\beta,\gamma\in\{1,\ldots,s\}$.
\begin{enumerate}
	\item $R(X,Y,\xi_\alpha,Z)=\varepsilon_\alpha\left\{\bar\eta(X)g(\varphi Y,\varphi Z)-\bar\eta(Y)g(\varphi X,\varphi Z)\right\}$;
	\item $R(\xi_\beta,Y,\xi_\alpha,Z)=\varepsilon_\beta\varepsilon_\alpha g(\varphi Y,\varphi Z)$;
	\item $R(\xi_\beta,\xi_\gamma,\xi_\alpha,Z)=0$;
	\item $R(\varphi X,\varphi Y,\xi_\alpha,Z)=0$;
	\item $R(U,Y,V,Z)=\bar\eta(U)\bar\eta(V)g(\varphi Y,\varphi Z)$.
\end{enumerate}
where, $\varepsilon_\alpha=g(\xi_\alpha,\xi_\alpha)=\pm1$ for any $\alpha\in\{1,\ldots,s\}$, and $\bar\eta=\sum_{\alpha=1}^{s}\varepsilon_\alpha\eta^\alpha$.
\end{proposition}
\proof With straightforward calculations, using (\ref{eq:04}), one gets $(1)$. The identities $(2)$, $(3)$ and $(4)$ are easy consequences of $(1)$, while $(5)$ follows from $(2)$. \endproof

As a consequence of the above properties, computing the Ricci tensor field $Ric(X,\xi_\alpha)$, for any $X\in \Gamma(TM)$ and any $\alpha\in\{1,\ldots,s\}$, we get
\begin{equation}\label{RicciXi}
\begin{split}
\Ric(X,\xi_\alpha)&=\sum_{i=1}^n\varepsilon_i \{R(X,E_i,\xi_\alpha,E_i)+R(X,\varphi E_i,\xi_\alpha,\varphi E_i)\}\\
                 &\quad+\sum_{\beta=1}^s \varepsilon_\beta R(X,\xi_\beta,\xi_\alpha,\xi_\beta)\\
                 &=\sum_{i=1}^n\varepsilon_i\varepsilon_\alpha \bar\eta(X)\{g(\varphi E_i,\varphi E_i)+g(E_i,E_i)\}=2n\varepsilon_\alpha\bar\eta(X)
\end{split}
\end{equation}
where $(E_i,\varphi E_i,\xi_\beta)$, $i\in\{1,\ldots,n\}$ and $\beta\in\{1,\ldots,s\}$, is any local orthonormal $\varphi$-adapted frame. Hence, using argumentations similar to those in \cite{KT}, we see at once that indefinite $\mathcal{S}$-manifolds can not be Einstein. Therefore, we introduce the notion of $\eta$-Einstein condition on an indefinite $g.f.f$-manifold as follows.

\begin{definition}
An indefinite $g.f.f$-manifold $(M,\varphi,\xi_\alpha,\eta^\alpha,g)$ is said to be \emph{$\eta$-Einstein} if there exist two functions $h,k\in\mathfrak{F}(M)$ such that
\begin{equation}\label{eq:05}
\Ric(X,Y)=hg(\varphi X,\varphi Y)+ k \bar\eta(X)\bar\eta(Y),
\end{equation}
for any $X,Y\in \Gamma(TM)$, where $\bar\eta=\sum_{\alpha=1}^s\varepsilon_\alpha\eta^\alpha$.
\end{definition}

\begin{remark}
Using (\ref{RicciXi}), from (\ref{eq:05}) one deduces that a $(2n+s)$-dimensional indefinite $\mathcal{S}$-manifold $(M,\varphi,\xi_\alpha,\eta^\alpha,g)$ is $\eta$-Einstein if and only if (\ref{eq:05}) holds with $k=2n$, that is
\begin{equation}\label{ricciS}
\Ric(X,Y)=hg(\varphi X,\varphi Y)+ 2n\bar\eta(X)\bar\eta(Y).
\end{equation}
for any $X,Y\in \Gamma(TM)$.
\end{remark}

\begin{remark}
It is easy to see that our definition reduces to the one given in \cite{KT}, when the signature of the metric is Euclidean. Indeed, in this case, we have $\varepsilon_\alpha=1$, for any $\alpha\in\{1,\ldots,s\}$, and (\ref{ricciS}) perfectly agrees with the condition $(1.12)$ of \cite{KT}, up to a multiplying factor.

Nevertheless, it is worth noting that (\ref{ricciS}) cannot be obtained from $(1.12)$ of \cite{KT} simply by inserting the $\varepsilon_\alpha$'s. Indeed, referring to \cite{KT} where the authors denote by $\tilde S$ the Ricci tensor field, if we replace each $\tilde\eta_x$ by $\varepsilon_x\tilde\eta_x$ in $(1.12)$, then we will get
\begin{equation*}
\begin{split}
\tilde S(\tilde X,\tilde Y)&=a(\tilde G(\tilde X,\tilde Y)+\sum_{x\neq y}\varepsilon_x\tilde\eta_x(\tilde X)\varepsilon_y\tilde\eta_y(\tilde X))\\
                           &+b(\sum_x\tilde\eta_x(\tilde X)\tilde\eta_x(\tilde Y)+\sum_{x\neq y}\varepsilon_x\tilde\eta_x(\tilde X)\varepsilon_y\tilde\eta_y(\tilde X)),
\end{split}
\end{equation*}
with $a+b=2n$ (up to a multiplying factor). The above expression is not equivalent to (\ref{ricciS}), due to the term $\sum_x\tilde\eta_x(\tilde X)\tilde\eta_x(\tilde Y)$, which does not agree with the analogous term obtained from \eqref{ricciS} by expanding it with the use of \eqref{eq:03}.
\end{remark}

\begin{remark}\label{tau}
Let $(M,\varphi,\xi_\alpha,\eta^\alpha,g)$ be an $\eta$-Einstein indefinite $\mathcal{S}$-manifold. Then the scalar curvature $\tau$ is given by
\begin{equation*}
\begin{split}
\tau&=\sum_{i=1}^n\varepsilon_i\{\Ric(E_i,E_i)+\Ric(\varphi E_i,\varphi E_i)\}+\sum_{\beta=1}^s\varepsilon_\beta \Ric(\xi_\beta,\xi_\beta)\\
    &=2nh+2n\sum_{\beta=1}^s\varepsilon_\beta=2n(h+\varepsilon),
\end{split}
\end{equation*}
where $(E_i,\varphi E_i,\xi_\beta)$, $i\in\{1,\ldots,n\}$ and $\beta\in\{1,\ldots,s\}$, is any local orthonormal $\varphi$-adapted frame and $\varepsilon=\sum_{\beta=1}^s\varepsilon_\beta$.
\end{remark}

Now we state the first Schur-type lemma for an $\eta$-Einstein indefinite $\mathcal{S}$-manifold. 

\begin{theorem}\label{T1}
Let $(M^{2n+s},\varphi,\xi_\alpha,\eta^\alpha,g)$, $n\geqslant2$ and $s\geqslant1$, be an $\eta$-Einstein indefinite $\mathcal{S}$-manifold. Then the function $h$ in (\ref{ricciS}) is constant.
\end{theorem}
\proof
Given $p\in M$, let $(E_i)_{i\in \{1,\ldots,2n+s\}}$ be a local orthonormal frame on a neighborhood $\mathfrak{U}$ of $p$ such that $(\nabla_{E_i}E_j)_p=0$, for any $i,j\in \{1,\ldots,2n+s\}$. Then, the Second Bianchi Identity, evaluated at the point $p$, has the form $\sigma_{(m,i,j)}E_m(R(E_i,E_j,E_k,E_l))=0$, for any $m,i,j,k\in\{1,\ldots,2n+s\}$. Putting $i=k$ and $j=l$, multiplying by $\varepsilon_i=g(E_i,E_i)$ and taking the sum over all $i\in\{1,\ldots,2n+s\}$, we get $E_m(\Ric(E_j,E_j)-2E_j(\Ric(E_j,E_m))=0$, for any $m,j\in\{1,\ldots,2n+s\}$. Multiplying again by $\varepsilon_j$ and taking the sum over all $j\in\{1,\ldots,2n+s\}$, by Remark \ref{tau}, we obtain
\begin{equation}\label{diffeq}
2nE_m(h)-2\sum_{j=1}^{2n+s}\varepsilon_jE_j(\Ric(E_j,E_m))=0.
\end{equation}
On the other hand, by \eqref{ricciS}, one has
\begin{equation}\label{eq:07}
\begin{split}
E_j(\Ric(E_j,E_m))&=E_j(h)g(\varphi E_j,\varphi E_m)+ h E_j(g(\varphi E_j,\varphi E_m))\\
                  &\quad+2n E_j(\bar\eta(E_j)\bar\eta(E_m)).
\end{split}
\end{equation}
Let us now calculate each term of the above identity separately. About the first one, using \eqref{eq:03}, we get, for any $m\in\{1,\ldots,2n+s\}$,
\begin{equation}\label{eq:08}
\sum_{j=1}^{2n+s}\varepsilon_jE_j(h)g(\varphi E_j, \varphi E_m)=E_m(h)-\sum_{\alpha=1}^s\eta^\alpha(E_m)\xi_\alpha(h).
\end{equation}
About the second term, by \eqref{eq:11} and \eqref{eq:12}, we have, at the point $p$, $E_j(\eta^\alpha(E_j))=-g(E_j,\varphi E_j)=0$. Using \eqref{eq:03} again, we get, for any $m\in\{1,\ldots,2n+s\}$,
\begin{equation}\label{eq:09}
\begin{split}
\sum_{j=1}^{2n+s}\varepsilon_jE_j(g(\varphi E_j,\varphi E_m))&=
-\sum_{j=1}^{2n+s}\sum_{\alpha=1}^{s}\varepsilon_j\varepsilon_\alpha E_j(\eta^\alpha(E_j)\eta^\alpha(E_m))\\
  &=-\sum_{j=1}^{2n+s}\sum_{\alpha=1}^{s}\varepsilon_jg(E_j,\xi_\alpha)g(\varphi E_m,E_j)\\
  &=-g(\bar\xi,\varphi E_m)=0.
\end{split}
\end{equation}
About the third term, since $E_j(\bar\eta(E_j))=0$, we have, for any $m\in\{1,\ldots,2n+s\}$,
\begin{equation}\label{eq:10}
\begin{split}
\sum_{j=1}^{2n+s}\varepsilon_jE_j(\bar\eta(E_j)\bar\eta(E_m))&=
\sum_{j=1}^{2n+s}\sum_{\alpha=1}^{s}\varepsilon_\alpha\varepsilon_jg(E_j,\bar\xi)g(\varphi E_m, E_j)\\
  &=\varepsilon g(\bar\xi,\varphi E_m)=0.
\end{split}
\end{equation}
Therefore, by \eqref{eq:07}, \eqref{eq:08}, \eqref{eq:09} and \eqref{eq:10}, \eqref{diffeq} yields
\[
(n-1)E_m(h)+\sum_{\alpha=1}^s\xi_\alpha(h)\eta^\alpha(E_m)=0,
\]
for any $m\in\{1,\ldots,2n+s\}$, from which it follows $(n-1)dh+\sum_{\alpha=1}^s\xi_\alpha(h)\eta^\alpha=0$. Applying this $1$-form to $\xi_\beta$, $\beta\in\{1,\ldots,s\}$, one obtains $\xi_\beta(h)=0$, for any $\beta\in \{1,\ldots,s\}$. Hence, since $n\geqslant2$ and $(n-1)X(h)=0$ for any $X\in\operatorname{Im}\varphi$, the claim follows. \endproof

\section{Indefinite $\mathcal{S}$-space forms as $\eta$-Einstein manifolds}
Let $(M,\varphi,\xi_\alpha,\eta^\alpha,g)$ be an indefinite $\mathcal{S}$-manifold and $p\in M$. For any non-lightlike unit vector $x\in\operatorname{Im}(\varphi_p)$, the sectional curvature of $(M,g)$ at $p$ with respect to the plane $\pi=span\{x,\varphi x\}$ is called, by definition, the \emph{$\varphi$-sectional curvature of $M$ at $p$, with respect to the $\varphi$-plane $\pi$}. When it is independent of the choice of the $\varphi$-plane at any point, the manifold $M$ is said to have \emph{pointwise constant $\varphi$-sectional curvature}. An indefinite $\mathcal{S}$-manifold with pointwise constant $\varphi$-sectional curvature is said to be an \emph{indefinite $\mathcal{S}$-space form} if the $\varphi$-sectional curvature does not depend on the point.

In \cite{LP} it is shown that an indefinite $\mathcal{S}$-manifold has pointwise constant $\varphi$-sectional curvature $c\in\mathfrak{F}(M)$ if, and only if, the Riemannian $(0,4)$-type curvature tensor field $R$ of $M$ satisfies the following identity
\begin{equation}\label{R}
\begin{split}
	R(X,Y,Z,W)&=\frac{c+3\varepsilon}{4}\{g(\varphi X,\varphi Z)g(\varphi Y,\varphi W)-g(\varphi Y,\varphi Z)g(\varphi X,\varphi W)\}\\
	&\quad\quad+\frac{c+\varepsilon}{4}\left\{\Phi(X,Z)\Phi(Y,W)-\Phi(Y,Z)\Phi(X,W)\right.\\
	&\qquad\qquad\left.+2\Phi(X,Y)\Phi(Z,W)\right\}\\
	&\quad\quad+\{\bar\eta(X)\bar\eta(Z)g(\varphi Y,\varphi W)-\bar\eta(Y)\bar\eta(Z)g(\varphi X,\varphi W)\\
	&\quad\quad+\bar\eta(Y)\bar\eta(W)g(\varphi X,\varphi Z)-\bar\eta(X)\bar\eta(W)g(\varphi Y,\varphi Z)\},
\end{split}
\end{equation}
for any $X,Y,Z,W\in\Gamma(TM)$.

Many examples of $\mathcal{S}$-manifolds with indefinite metrics have been introduced and studied in different contexts. Namely, it is possible to endow $\mathbb{R}^4$, $\mathbb{R}^6$ and $U(2)$ with non-trivial indefinite $\mathcal{S}$-structures. In particular, $\mathbb{R}^4$ and $U(2)$ turn out to be both Lorentzian $\mathcal{S}$-space forms, and it is easy to check that they are $\eta$-Einstein with $h=0$ and $h=4$, respectively (see \cite{LP} for more details about the non-compact examples, and \cite{LP2} for the $U(2)$ case).

Now, we are going to show that any indefinite $\mathcal{S}$-manifold with pointwise constant $\varphi$-sectional curvature is $\eta$-Einstein.

\begin{theorem}
Let $(M^{2n+s},\varphi,\xi_\alpha,\eta^\alpha,g)$, $n\geqslant2$ and $s\geqslant1$, be an indefinite $\mathcal{S}$-manifold with pointwise constant $\varphi$-sectional curvature $c\in\mathfrak{F}(M)$. Then $M$ is $\eta$-Einstein.
\end{theorem}

\proof
Let $(E_i,\varphi E_i,\xi_\alpha)$, $i\in\{1,\ldots,n\}$ and $\alpha\in\{1,\ldots,s\}$, be a local orthonormal $\varphi$-adapted frame. We have
\[
\Ric(X,Y)=\sum_{i=1}^n\varepsilon_i \{R(X,E_i,Y,E_i)+R(X,\varphi E_i,Y,\varphi E_i)\}+\sum_{\alpha=1}^s\varepsilon_\alpha R(X,\xi_\alpha,Y,\xi_\alpha).
\]
for any $X,Y\in\Gamma(TM)$. Using (\ref{R}) we have, for any $i\in\{1,\ldots,n\}$ and any $\alpha\in\{1,\ldots,s\}$:
\begin{align*}
	R(X,E_i,Y,E_i)&=\frac{c+3\varepsilon}{4}\,\{g(\varphi X,\varphi Y)\varepsilon_i-g(\varphi E_i,\varphi Y)g(\varphi X,\varphi E_i)\}\\
	&\quad\quad +3\,\frac{c-\varepsilon}{4}\,\Phi(X,E_i)\Phi(Y,E_i)+\bar\eta(X)\bar\eta(Y)\varepsilon_i,\\
	R(X,\varphi E_i,Y,\varphi E_i)&=\frac{c+3\varepsilon}{4}\,\{g(\varphi X,\varphi Y)\varepsilon_i-g(E_i,\varphi Y)g(\varphi X, E_i\}\\
	&\quad\quad +3\,\frac{c-\varepsilon}{4}\,\Phi(X,\varphi E_i)\Phi(Y,\varphi E_i)+\bar\eta(X)\bar\eta(Y)\varepsilon_i,\\
	R(X,\xi_\alpha,Y,\xi_\alpha)&=g(\varphi X,\varphi Y).
\end{align*}
Therefore, for any $X,Y\in\Gamma(TM)$, we get
\begin{align}\label{Ricci01}
\Ric(X,Y)&=2n\,\frac{c+3\varepsilon}{4}\,g(\varphi X,\varphi Y)\nonumber\\
        &\quad\quad-\frac{c+3\varepsilon}{4}\sum_{i=1}^n\varepsilon_i\{g(\varphi E_i,\varphi Y)g(\varphi X,\varphi E_i)
                                                           +g(\varphi X,\varphi E_i)g(\varphi Y,\varphi E_i)\}\nonumber\\
        &\quad\quad+3\,\frac{c-\varepsilon}{4}\sum_{i=1}^n\varepsilon_i\{g(\varphi X,E_i)g(\varphi Y,E_i)
                                                                  +g(\varphi X,\varphi E_i)g(\varphi Y,\varphi E_i)\}\nonumber\\
        &\quad\quad+2n\bar\eta(X)\bar\eta(Y)+\varepsilon g(\varphi X,\varphi Y)\\
         &=\frac{c+3\varepsilon}{4}\,(2n-1)g(\varphi X,\varphi Y)+3\,\frac{c-\varepsilon}{4}\,g(\varphi X,\varphi Y)\nonumber\\
	      &\quad\quad+2n\bar\eta(X)\bar\eta(Y)+\varepsilon g(\varphi X,\varphi Y)\nonumber\\
	       &=\frac{1}{2}\,\{n(c+3\varepsilon)+c-\varepsilon\}g(\varphi X,\varphi Y)+2n\bar\eta(X)\bar\eta(Y).\nonumber
\end{align}
Then $M$ turns out to be $\eta$-Einstein. \endproof

Clearly, for an indefinite $\mathcal{S}$-manifold with pointwise constant $\varphi$-sectional curvature $c\in\mathfrak{F}(M)$, \eqref{Ricci01} yields \eqref{ricciS} with $h=\frac{1}{2}\,\{n(c+3\varepsilon)+c-\varepsilon\}$ and Theorem \ref{T1} implies the following consequence.

\begin{theorem}
Let $(M^{2n+s},\varphi,\xi_\alpha,\eta^\alpha,g)$, $n\geqslant2$ and $s\geqslant1$, be an indefinite $\mathcal{S}$-manifold with pointwise constant $\varphi$-sectional curvature $c\in\mathfrak{F}(M)$. Then $c$ is a constant function on $M$, i.e. $M$ is an indefinite $\mathcal{S}$-space form. 
\end{theorem}

The above result extends the ones of \cite{KT} to the semi-Riemannian setting. To conclude we want to give the following remark that states a relation between the $\eta$-Einstein notion on an indefinite $\mathcal{S}$-manifold and the K\"ahler-Einstein one.
\begin{remark}
In \cite{LP2} it is stated that an indefinite K\"ahler structure $(J,g')$ on a manifold $N$ can be lifted to an indefinite $\mathcal{S}$-structure $(\varphi,\xi_\alpha,\eta^\alpha,g)$ on the total space $M$ of a principal toroidal bundle, whose projection $\pi:(M,\varphi,\xi_\alpha,\eta^\alpha,g)\rightarrow(N,J,g')$ turns out to be a semi-Riemannian submersion with totally geodesic fibres. Looking at \cite[p.\ 15 and p.\ 145]{FIP}, in this context the Ricci formulas yield
\begin{equation*}
	\Ric(X,Y)=\Ric'(X',Y')\circ\pi-2g(\bar\xi,\bar\xi)g(\varphi X,\varphi Y);
\end{equation*}
where $X,Y$ are basic vector fields $\pi$-related to $X',Y'$. When $N$ is an Einstein manifold, by the above formula, it is clear that $M$ is an $\eta$-Einstein manifold.
\end{remark}

\end{document}